\documentclass[
	a4paper,
	12pt
	]{article}

\usepackage[utf8]{inputenc}
\usepackage{hyperref}
\usepackage{libertine}

\usepackage[dvipsnames]{xcolor}

\usepackage{geometry}
\usepackage{setspace}
\spacing{1.1}
\usepackage{caption}
\usepackage{soul}

\usepackage{graphicx}
\usepackage{enumerate}

\usepackage{amsmath,amsfonts,amssymb,xfrac,nicefrac}
\usepackage{amsthm}
\usepackage{mathrsfs}
\usepackage{braket}
\usepackage{aliascnt}
\usepackage{microtype}
\usepackage{mathtools}
\usepackage{etoolbox}
\usepackage{mathrsfs}

\usepackage{comment}

\usepackage{titlesec}
\renewcommand{\thesection}{\arabic{section}}
\renewcommand{\thesubsection}{\arabic{section}.\arabic{subsection}}
\titlelabel{\thetitle. }
\titleformat{\section}{\center\large\bfseries\boldmath}{\thesection.}{1ex}{}
\titleformat{\subsection}{\center\bfseries\boldmath}{\thesubsection.}{1ex}{}

\newcommand{\periodafter}[1]{\ifstrempty{#1}{}{#1.}}
\titleformat{\subsection}[runin]
{\bfseries}{\thesubsection.}{1ex}{\periodafter}

\usepackage{algorithm}
\usepackage{algpseudocode}
	\algrenewcommand\algorithmicrequire{\textbf{Input:}}
	\algrenewcommand\algorithmicensure{\textbf{Output:}}
	\newcommand{\Goto}{\textbf{goto}}
	
	\algblockdefx{ParallelFor}{EndParallelFor}%
		[1]{\textbf{for} #1 \textbf{parallel do}}%
		{\textbf{end parallel do}}%

		\errorcontextlines\maxdimen

		\makeatletter
		\newcommand*{\algrule}[1][\algorithmicindent]{\makebox[#1][l]{\hspace*{.5em}\vrule height .75\baselineskip depth .25\baselineskip}}%

		\newcount\ALG@printindent@tempcnta
		\def\ALG@printindent{%
		    \ifnum \theALG@nested>0% is there anything to print
		        \ifx\ALG@text\ALG@x@notext% is this an end group without any text?
		            \addvspace{-3pt}% FUDGE for cases where no text is shown, to make the rules line up
		        \else
		            \unskip
		            \ALG@printindent@tempcnta=1
		            \loop
		                \algrule[\csname ALG@ind@\the\ALG@printindent@tempcnta\endcsname]%
		                \advance \ALG@printindent@tempcnta 1
		            \ifnum \ALG@printindent@tempcnta<\numexpr\theALG@nested+1\relax% can't do <=, so add one to RHS and use < instead
		            \repeat
		        \fi
		    \fi
		    }%
		\usepackage{etoolbox}
		\patchcmd{\ALG@doentity}{\noindent\hskip\ALG@tlm}{\ALG@printindent}{}{\errmessage{failed to patch}}
		\makeatother
\usepackage{tikz}
\usepackage{tikz-cd}
	\usetikzlibrary{arrows,calc}
	\usetikzlibrary{decorations.pathreplacing,calc}
	\newcommand{\tikzmark}[1]{\tikz[overlay,remember picture] \node (#1) {};}

\usepackage{hyperref}
	\hypersetup{
	    colorlinks,
	    linkcolor=darkblue,
	    citecolor=darkblue,
	    urlcolor=darkblue}
	    \definecolor{darkblue}{rgb}{0.0,0.0,0.5}
	\definecolor{darkgreen}{rgb}{0.0,0.5,0.0}
	\definecolor{darkred}{rgb}{0.8,0.3,0.0}

\newcommand{\Date}{15~March~2022}
\newcommand{\NumberMSASequences}{8,297,154~}
\newcommand{\NumberMSACleanSequences}{5,323,639~}
\newcommand{\NumberMSADistinctSequences}{359,650~}
\newcommand{\LengthMSACleanSequences}{4,874nt}

\newcommand{\episetDOI}{\url{https://doi.org/10.55876/gis8.220629ug}}
\newcommand{\Alignment}{\texttt{msa\_0315.fasta}}
\newcommand{\DownloadDate}{28~March~2022}

\usepackage[capitalise]{cleveref}
\usepackage{enumerate}

\newtheorem*{theorem*}{Theorem}
\newtheorem{theorem}{Theorem}[section]

\theoremstyle{definition}

\makeatletter
\let\@fnsymbol\@arabic
\makeatother

\title{\vspace{-1.7cm} \Large\bfseries MuRiT: Efficient Computation of Pathwise Persistence Barcodes in Multi-Filtered Flag Complexes via \\ Vietoris-Rips Transformations}

\author{Maximilian Neumann \!\footnote{Mathematics Department, Karlsruhe Institute of Technology, Karlsruhe, Germany} \,, Michael Bleher \!\footnote{Mathematical Institute, Heidelberg University, Heidelberg, Germany} \,, Lukas Hahn\textsuperscript{\,2}, Samuel Braun \!\footnote{Steinbuch Centre for Computing, Karlsruhe Institute of Technology, Karlsruhe, Germany} \,, \\ Holger Obermaier\textsuperscript{\,3}, Mehmet Soysal\textsuperscript{\,3}, Ren\'e Caspart\textsuperscript{\,3}, Andreas Ott\textsuperscript{\,1,2}}

\begin{document}

\maketitle

\vspace{-5mm}

\begin{abstract}
\noindent Multi-parameter persistent homology naturally arises in applications of persistent topology to data that come with extra information depending on additional parameters, like for example time series data.
We introduce the concept of a Vietoris-Rips transformation, a method that reduces the computation of the one-parameter persistent homology of pathwise subcomplexes in multi-filtered flag complexes to the computation of the Vietoris-Rips persistent homology of certain semimetric spaces.
The corresponding pathwise persistence barcodes track persistence features of the ambient multi-filtered complex and can in particular be used to recover the rank invariant in multi-parameter persistent homology. We present \texttt{MuRiT}, a scalable algorithm that computes the pathwise persistence barcodes of multi-filtered flag complexes by means of Vietoris-Rips transformations. Moreover, we provide an efficient software implementation of the \texttt{MuRiT} algorithm which resorts to \texttt{Ripser} for the actual computation of Vietoris-Rips persistence barcodes. To demonstrate the applicability of \texttt{MuRiT} to real-world datasets, we establish \texttt{MuRiT} as part of our \texttt{CoVtRec} pipeline for the surveillance of the convergent evolution of the coronavirus SARS-CoV-2 in the current COVID-19 pandemic.
\end{abstract}

\section{Introduction}
\label{sec:intro}
Persistent homology is one of the most important tools in computational topology and topological data analysis.
It has the capability to detect and explore qualitative features of complex datasets that are encoded in the geometric shape of the dataset and are otherwise hard to extract with traditional methods (see e.g.~\cite{Edelsbrunner2008, Carlsson2009, Otter2017, Ghrist2007, Weinberger2011, EH2010, Chazal2016, Oudot2015, Dey2022}).
A common approach is Vietoris-Rips persistent homology, which analyzes the geometric shape of metric datasets at varying distance scales.
In many applications, however, data points come with extra information that is given in terms of additional attributes and one wishes to leverage this extra information in the topological data analysis.
A typical example of this is time series data.

Our motivating application in this paper is exactly of this sort---we use persistent homology for the surveillance of emerging adaptive mutations in the evolution of the coronavirus SARS-CoV-2 in the current COVID-19 pandemic \cite{topologyidentifies2022, covtrec}.
Here the dataset consists of coronavirus gene sequences.
The use of persistent homology to analyze the evolution of viruses was initiated by Chan, Carlsson and Rabad\'an \cite{Chan2013}.
The coronavirus adapts itself to the human host by developing new variants by mutating its genome.
In \cite{topologyidentifies2022} we introduced a topological descriptor
for the adaptiveness of a given mutation in the genome of the coronavirus that is defined by counting certain one-dimensional cycle representatives in the Vietoris-Rips persistent homology of the gene sequences dataset (see \cref{sec:application}).
Now each coronavirus gene sequence in the dataset is assigned the date at which it was collected from a patient.
In this way, the dataset comes with a natural stratification by sampling time, with a bunch of new sequences being added every day.

Ideally, one would like to exploit this additional information and monitor topological signals of adaptation over time in order to tell whether a given mutation is likely to become adaptive in the future.
A naive approach is to regard time as an external parameter, and to run the persistence analysis separately for each sub-dataset consisting of all sequences that have been collected up to a given point in time.
However, this approach is computationally expensive, as the whole analysis has to be repeated many times.
Moreover, classes in persistent homology computed at different time steps will in general not be related with each other.
As we will see, all these issues can be resolved by including time as an additional parameter into the persistence analysis itself.
The natural setup for this is multi-parameter persistent homology of multi-filtered simplicial complexes introduced by Carlsson and Zomorodian
\cite{Carlsson2009multidimensionalpersistence, Carlsson2009compute}. 
While it is a challenge to compute multi-parameter persistent homology in general \cite{introtomultipers}, it turns out that for our applications in viral evolution one only needs to compute the persistent homology of certain one-filtered subcomplexes in multi-filtered flag complexes.

In the present paper, we address this problem and present \texttt{MuRiT}, a fast and scalable algorithm for the computation of the persistent homology of arbitrary one-filtered subcomplexes of a given multi-filtered flag complex (see \cref{sec:MuRiT}).
The main idea of the \texttt{MuRiT} algorithm is to apply \emph{Vietoris-Rips transformations} in order to reduce the computation of the persistent homology of one-filtered subcomplexes in a multi-filtered flag complex to the computation of the usual Vietoris-Rips persistent homology of certain semimetric spaces.
We will explain Vietoris-Rips transformations in more detail in the next paragraph.
The actual computation of the Vietoris-Rips persistent homology of the semimetric space can then be carried out independently with basically any of the presently available software packages \cite{Otter2017}, depending on the needs of the particular application one has in mind. However, one has to make sure that the chosen software package is able to handle the Vietoris-Rips persistent homology of \emph{semimetric} spaces that do not necessarily satisfy the triangle inequality.

We provide an efficient software implementation of the \texttt{MuRiT} algorithm at \url{https://github.com/tdalife/murit}.
In its current form, this implementation is tailored to the case of Vietoris-Rips persistent homology of multi-filtered point cloud datasets, a setup which naturally arises in the Vietoris-Rips persistence analysis of time series data.
By default, our implementation of \texttt{MuRiT} resorts to the \texttt{Ripser} software package by Bauer \cite{bauer2021ripser} for the actual computation of persistence barcodes.
Note at this point that \texttt{Ripser} is able to compute the Vietoris-Rips persistence barcodes also for semimetric spaces that do not necessarily satisfy the triangle inequality \cite{notriangle}.
In this way, \texttt{MuRiT} takes full advantage of the computational power of \texttt{Ripser}, which is among the most efficient implementations for the computation of persistent homology to date \cite{Otter2017}.
\texttt{MuRiT} is part of our \texttt{CoVtRec} pipeline for the surveillance of potentially adaptive mutations in the evolution of the coronavirus SARS-CoV-2 in the current COVID-19 pandemic \cite{covtrec} (see \cref{sec:covtrec}).
Thanks to highly optimized algorithms that take advantage of the tree-like structure of the gene sequences dataset \cite{gromovhyper}, \texttt{CoVtRec} has the capability to process very large SARS-CoV-2 genomic datasets and easily scales to hundreds of thousands of distinct genomes.

\begin{figure}[t]
	\centering
	\includegraphics[width=0.7\textwidth]{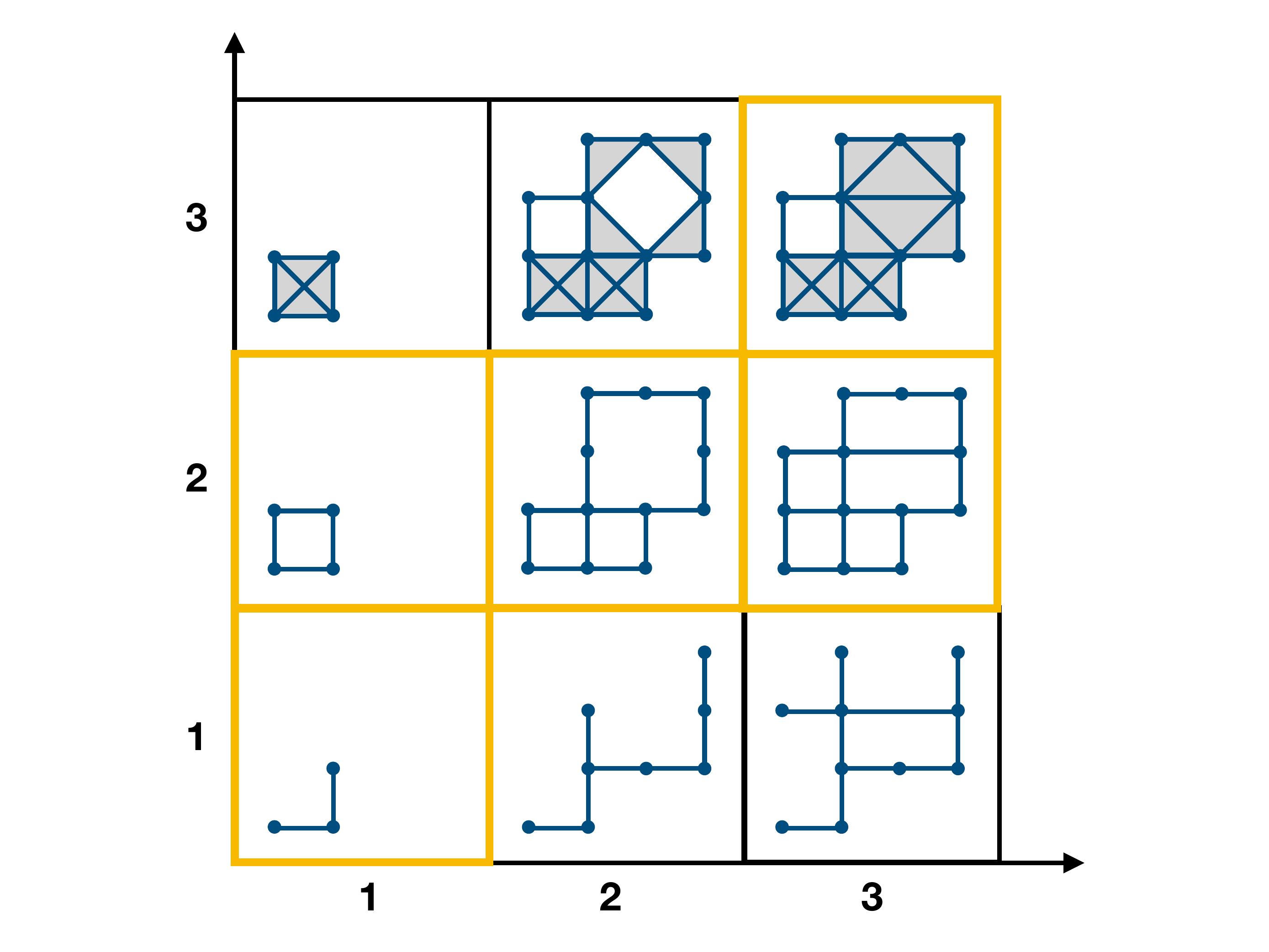}
\captionsetup{labelfont=bf}
\caption{\textbf{Example of a multi-filtered flag complex.}
The displayed flag complex $X$ is bi-filtered with three filtration steps in each dimension. The yellow squares mark the one-filtered subcomplex $X_{(1,1)} \subseteq X_{(1,2)} \subseteq X_{(2,2)}  \subseteq  X_{(3,2)} \subseteq X_{(3,3)}$ of $X$. }
\label{fig:1}
\end{figure}

Let us state our main result and outline the basic idea underlying Vietoris-Rips transformations (see \cref{sec:main-results}).
Assume that $X$ is a finite $P$-filtered flag complex for some partially ordered set $P=(P, \leq)$, and consider a \emph{(discrete) path} in $P$ that is given by a monotone sequence ${\nu=(\nu_1 \leq \nu_2 \leq \nu_3 \leq \dots  )}$ of elements in $P$. This gives rise to a one-filtered subcomplex ${X_{\nu}=(X_{\nu_1} \subseteq X_{\nu_2} \subseteq X_{\nu_3} \subseteq \dots)}$ of $X$ (see Figure~\ref{fig:1}). Ideally, for the actual computation of the one-parameter persistent homology of~$X_{\nu}$ we would like to resort to any of the currently available efficient algorithms for the computation of Vietoris-Rips persistent homology, like for example \texttt{Ripser}.
To that end, we construct a semimetric~$d$ on the vertex set $\text{Vert}(X_{\nu})$ that encodes the filtration steps of the one-filtered complex $X_{\nu}$ in a suitable way, and define the \emph{Vietoris-Rips transformation} of~$X_{\nu}$ as the Vietoris-Rips complex
\[
\widehat{\text{VR}} (X_{\nu}) := \text{VR}(\text{Vert}(X_{\nu}), d) )
\]
of the semimetric space $(\text{Vert}(X_{\nu}), d)$.
Then we prove that the one-parameter persistent homology of the Vietoris-Rips transformation $\widehat{\text{VR}} (X_{\nu})$ recovers the persistent homology of the original filtration~$X_{\nu}$ in the sense that there is an isomorphism
\[
H_{\ell} (X_{\nu}) \cong H_{\ell} \big (  \widehat{\text{\normalfont VR}}(X_{\nu}) \big)
\]
of persistence modules in every positive degree $\ell > 0$ (see \cref{mainthm}).
In the proof we use the fact that the subcomplexes $X_{\nu_i}$ are flag complexes.
Let us remark that the distance function $d$ will in general not satisfy the triangle inequality, which is why in the definition of the Vietoris-Rips transformation we need to consider Vietoris-Rips complexes of semimetric spaces.
We will normally phrase our result in terms of \emph{pathwise barcodes} by saying that in every positive degree~$\ell > 0$, the persistence barcode of the complex $X$ along the path $\nu$ is computed by
\[
\mathcal{B}_{  \ell}  (X_{\nu}) = \mathcal{B}_{ \ell}  \big (  \widehat{\text{\normalfont VR}} (X_{\nu}   ) \big )
\]
in terms of the usual Vietoris-Rips persistence barcode of the Vietoris-Rips transformation $\widehat{\text{\normalfont VR}} (X_{\nu})$ of the one-filtered subcomplex $X_{\nu}\subseteq X$ (see \cref{sec:PathwiseBarcodes}). 

As is shown in \cite{Carlsson2009multidimensionalpersistence}, there exists no discrete and complete invariant in multi-parameter persistent homology like the persistence barcode known from one-parameter persistence. But there are several approaches to define invariants for multi-persistence, like for example the \emph{rank invariant} introduced in \cite{Carlsson2009multidimensionalpersistence}, \textit{Hilbert functions} (see e.g.~\cite{introtomultipers}), \textit{multi-graded Betti numbers} (see e.g.~\cite{Sturmfels2005}), \textit{signed barcodes} \cite{signedbarcodes} and \textit{fibered barcodes} \cite{lesnick2015interactive, introtomultipers,Bettinumbersmultipers}.
Our approach to consider pathwise barcodes of finite multi-filtered flag complexes is reminiscent of the concept of fibered barcodes, where the basic idea is to compute persistence barcodes along affine lines in $\mathbb{R}^n$.
In particular, both pathwise and fibered barcodes recover the rank invariant for multi-parameter persistence (see \cref{subsec:RankInvariant}).
With \texttt{RIVET}, Wright, Lesnick et al. \cite{RIVET,lesnick2015interactive} provide a software package for working with two-parameter persistent homology, which provides an interactive visualization of the Hilbert function, the bi-graded Betti numbers, and the fibered barcode.
Another algorithm specifically designed for the efficient computation of the persistent homology of directed flag complexes is the \texttt{Flagser} software package by L\"utgehetmann, Govc, Smith and Levi \cite{Flagser}.
A particular feature of \texttt{MuRiT} in comparison with \texttt{Flagser} is that it does not do the actual computation of persistent homology by itself.
In this way, \texttt{MuRiT} offers maximum flexibility regarding the choice of software package for the computation of persistent homology.
This feature is, for example, indispensable in our application of \texttt{MuRiT} to the evolution of the coronavirus, as we need to use a custom version of \texttt{Ripser} that is specifically optimized for the efficient localization of cycles in the gene sequences dataset.
Another feature of our implementation of \texttt{MuRiT} is that it is genuinely designed to deal with pathwise subfiltrations of multi-filtered complexes.

The paper is organized as follows. In \cref{sec:preliminaries}, we fix the notation and recall some basic facts and definitions about the persistent homology of multi-filtered flag complexes.
In \cref{sec:main-results}, we introduce the notion of a Vietoris-Rips transformation, define pathwise barcodes and relate them with the rank invariant, and  present the \texttt{MuRiT} algorithm.
The final \cref{sec:application} discusses an application of the \texttt{MuRiT} algorithm to the evolution of the coronavirus.

\subsection*{Acknowledgements} The authors gratefully acknowledge all data contributors, i.e.~the Authors and their Originating laboratories responsible for obtaining the specimens, and their Submitting laboratories for generating the genetic sequence and metadata and sharing via the GISAID Initiative \cite{Shu2017, Khare2021}, on which this research is based.
An acknowledgement table is accessible online at \episetDOI.
The authors acknowledge the use of de.NBI Cloud and the support by the High Performance and Cloud Computing Group at the Zentrum für Datenverarbeitung of the University of Tübingen and the German Federal Ministry of Education and Research (BMBF) through grant no 031 A535A. They thank M. Hanussek for IT support and early access to VALET \cite{valet}. The authors further acknowledge support from the Interdisciplinary Center for Scientific Computing at Heidelberg University and the development work of the Scientific Software Center of Heidelberg University carried out by L. Keegan and D. Kempf \cite{hammingdist}. A.O. acknowledges funding by the Federal Ministry of Education and Research (BMBF) and the Baden-Württemberg Ministry of Science as part of the Excellence Strategy of the German Federal and State Governments (KIT~Centers, "Topological Genomics"). A.O. and M.N. acknowledge funding by the Vector Foundation ("Topological Genomics"). L.H. and M.B. were supported by the Deutsche Forschungsgemeinschaft (DFG, German Research Foundation) under Germany’s Excellence Strategy EXC 2181/1 - 390900948 (the Heidelberg STRUCTURES Excellence Cluster). L.H. thanks the Evangelisches Studienwerk Villigst for their support.

\subsection*{Author Contributions} M.N. developed the concept of Vietoris-Rips tranformations, which grew out of the applications part of M.N.'s Master's thesis under the supervision of A.O.; M.N., M.B. designed and developed the \texttt{MuRiT} algorithm; M.B. designed and implemented the software package \texttt{MuRiT}; M.N., M.B., L.H., A.O. designed and implemented \texttt{CoVtRec}; M.B., L.H., A.O. curated data for \texttt{CoVtRec}; M.N., M.B., L.H., A.O. performed computational analyses; M.N.,M.B., L.H.,  A.O., S.B., H.O., M.S., R.C. developed and implemented software for \texttt{CoVtRec}; M.B., L.H., A.O. acquired computing resources for \texttt{CoVtRec}; M.N., M.B., L.H., A.O. drafted the manuscript; all authors contributed to the final version of this article.
\vspace{1cm}
\section{Preliminaries}
\label{sec:preliminaries}
\subsection{Partially Ordered Sets}

Let us denote by $\mathbb{N}=\{1,2,3,\dots \}$ the set of natural numbers, and write $\mathbb{N}_{0}=\mathbb{N} \cup \{0\}$. For~$n \in \mathbb{N}$, we will be working with the following partial order on the $n$-fold cartesian product $\mathbb{R}^n$ that is induced by the usual total order $\leq$ on the real line $\mathbb{R}$. For any pair of tuples ${{a}=(a_1, \dots, a_n)}$ and ${b}=(b_1, \dots, b_n)$ in $\mathbb{R}^n$, we define ${a} \leq {b}$ if $a_i \leq b_i$ for all $i \in \{1, \dots, n\}$. The subset $\mathbb{N}^n \subseteq \mathbb{R}^n$ naturally becomes a partially ordered set with the induced partial order $\leq$ inherited from $\mathbb{R}^n$.

A poset $P=(P, \leq)$ is said to have \textit{dimension} $n$ if there exists an order preserving embedding $(P, \leq) \hookrightarrow (\mathbb{R}^m, \leq)$ for $m=n$, but at the same time no such embedding exists for $m<n$.

\subsection{Simplicial Complexes and Graphs}

We briefly recall some basic facts and definitions about simplicial complexes, and fix some notation and terminology.

An undirected \emph{graph} is a pair $G = (V, E)$ consisting of a set $V$ of \emph{vertices} and a set $E$ of unordered pairs of vertices in $V$ called the \emph{edges}.

An \emph{(abstract) simplicial complex} is a set $X$ of nonempty finite sets such that if $\sigma$ is an element of~$X$, so is every nonempty subset of $\sigma$.
The elements of $X$ are called the \emph{simplices} of the complex~$X$.
A simplex with $k+1$ elements is called a \emph{$k$-simplex}, and $k$ is also called its \textit{dimension}.
As a particular case of this, $0$-simplices in $X$ are also called \emph{vertices} and $1$-simplices in $X$ are called  \emph{edges}. Any subset of $X$ that is itself a simplicial complex is called a \emph{subcomplex} of $X$.
For every non-negative integer~$k$, the subcomplex $X^{(k)} \subseteq X$ consisting of all simplices in $X$ of dimension at most $k$ is called the \emph{$k$-skeleton} of $X$.
The $0$-skeleton $X^{(0)}$ is also called the \emph{vertex set} $\text{Vert}(X)$ of $X$, and the complement $\text{Edge}(X) := X^{(1)} \setminus X^{(0)}$ of the vertex set in the $1$-skeleton will be called the \emph{edge set} of $X$.

Consider an undirected graph $G = (V, E)$. A \emph{clique} in the graph $G$ is a finite subset $C \subseteq V$ of the set of vertices such that for any two distinct vertices $u$ and $v$ in $C$, the unordered pair $\{u, v\}$ formed by these two vertices is contained as an edge in $E$.
We denote by $\mathscr{C}(G)$ the set of all cliques in~$G$.
By construction, $\mathscr{C}(G)$ is a simplicial complex and is called the \emph{clique complex} of the graph $G$.

Let $X$ be a simplicial complex. Observe that the vertex and edge sets of $X$ give rise to an undirected graph \begin{equation}
G(X) := ( \text{Vert}(X), \text{Edge}(X) )  \notag
\end{equation}
The complex $X$ is called a \emph{flag complex} if it satisfies the condition $X = \mathscr{C}(G(X)).$
In other words, a flag complex is by definition the clique complex of the graph formed by its vertex and edge sets.
Then the simplices of the flag complex are precisely the cliques in its vertex set.
Note that in this way, every flag complex is completely determined by its $1$-skeleton.

\subsection{Filtered Sets and Filtered Simplicial Complexes} Let $P=(P, \leq)$ be a poset and $X$ be a set. A $P$-\emph{filtration} of~$X$ is a family of sets $X_{\bullet}=(X_p)_{p \in P}$ satisfying the following conditions:
\begin{enumerate}[(i)]
	\item $X_p \subseteq X$ is a subset for every $ p \in P$.
	\item $X_p \subseteq X_q$ for all $p,q \in P$ with $p \leq q$.
	\item $\bigcup_{p \in P} X_p=X.$
\end{enumerate}
A set $X$ is called $P$-filtered if it admits a $P$-filtration.
If $P$ is $n$-dimensional, $X$ is called \emph{$n$-filtered}.
$X$ is called \emph{multi-filtered} if $X$ is $n$-filtered for some $n\geq 2$.

A simplicial complex $X$ is called \emph{$P$-filtered} if $X$ is a $P$-filtered set and $X_p\subseteq X$ is a subcomplex for every $p \in P$. A $P$-filtered simplicial complex $X$ is called a \emph{$P$-filtered flag complex} if $X$ is a flag complex and the subcomplexes $X_p$ are flag complexes for all $p \in P$.

\subsection{Vietoris-Rips Complexes of Semimetric Spaces}

In the literature, Vietoris-Rips complexes are normally defined for metric spaces.
It is key to our approach in this paper to consider Vietoris-Rips complexes for a larger class of spaces equipped with a more general notion of distance function that is not required to satisfy the triangle inequality. Let $S$ be a set and $[0, \infty]=\mathbb{R}_{\geq 0} \cup \{\infty\}$, where $\infty$ henceforth denotes $+\infty$ for short.
A function $d\!:S \times  S \to [0, \infty]$ is called a \emph{semimetric} if it satisfies the following two axioms:
\begin{enumerate}[(i)]
	\item $d(x,y)=d(y,x)$ for all $x,y\in S$.
	\item $d(x,y)=0$ if and only if $x=y$, for all $x,y\in S$.
\end{enumerate}
The pair $(S,d)$ is called a \emph{semimetric space}. The function $d$ will also be called the \emph{distance function} of the semimetric space $(S,d)$. Note that a semimetric space is not required to satisfy the triangle inequality, and that we allow the distance function to take the value $\infty$.

Let $(S,d)$ be a semimetric space.
For every $r\in [0, \infty)$, the \emph{Vietoris-Rips complex of $(S,d)$ at scale $r$} is the abstract simplicial complex defined by
\[
\text{VR}_r(S,d) := \big\{ \sigma \subseteq S \,\big|\, \emptyset \neq \sigma \text{ finite and } d(x,y) \leq r \, \text{ for all } \,  x, y \in \sigma \big\}.
\]
Let us remark that this definition makes sense also if the distance function $d$ is not required to satisfy the triangle inequality.
Note moreover that $\text{VR}_r(S,d)$ is in fact a flag complex.
Its simplices are precisely all finite non-empty subsets of the set $S$ whose points have pairwise distance at most $r$.
We will also consider the simplicial complex
\begin{equation}
\text{VR}(S,d) := \{ \sigma \subseteq \text{Vert}(X) \,\big|\, \emptyset \neq \sigma \text{ finite and } d(x,y) < \infty \, \text{ for all } \,  x, y \in \sigma \big\}. \notag
\end{equation}
It is a flag complex and we will refer to it as the \emph{Vietoris-Rips complex} of the semimetric space~$(S,d)$.
Note that $\text{VR}(S,d)$ becomes a $[0,\infty)$-filtered flag complex with Vietoris-Rips filtration $\text{VR}_{\bullet}(S,d)=(\text{VR}_r(S,d) )_{r \in [0,\infty)}$.

\subsection{Persistence Modules}

Let $P=(P, \leq)$ be a poset, and fix a coefficient field $\mathbb{F}$. In later computations we will normally choose $\mathbb{F}=\mathbb{F}_p$ to be a finite field of prime order.
A \emph{persistence module} over $P$ is a functor
\[
M: P \to \mathbf{Vec}_{\mathbb{F}}, \quad p \mapsto M(p)
\]
from $P$ into the category $\mathbf{Vec}_{\mathbb{F}}$ of vector spaces over the field $\mathbb{F}$. It assigns to every pair $p,q \in P$ with $p \leq q$ an $\mathbb{F}$-linear map denoted by \begin{equation}
M(p \leq q): M(p) \to M(q) \notag.
\end{equation}
If $P$ is $n$-dimensional, then $M$ is also called an \emph{$n$-parameter persistence module}. $M$ is called a \emph{multi-parameter persistence module} if $M$ is an $n$-parameter persistence module for some $n\geq 2$. A persistence module $M$ over $P$ is called \textit{pointwise finite dimensional} if $M(p)$ is a finite dimensional $\mathbb{F}$-vector space for all $p \in P$.

Let now $X$ be a finite $P$-filtered simplicial complex. Then for every non-negative integer $\ell \geq 0$, the assignment
\[
H_{\ell}(X): P \to  \mathbf{Vec}_{\mathbb{F}}, \quad p \mapsto H_{\ell}(X_p)
\]
defines a pointwise finite dimensional persistence module over $P$, where $H_{\ell}(X_p)=H_{\ell}(X_p; \mathbb{F})$ denotes the $\ell$-th simplicial homology of $X_p$ with coefficients in the field $\mathbb{F}$. 

\section{Main Results}
\label{sec:main-results}

\subsection{Vietoris-Rips Transformations}
\label{sec:VRTransformations}
We start with a construction that assigns a finite semimetric space to any finite $\mathbb{N}$-filtered flag complex $X$. For this, we turn the vertex set of $X$ into a semimetric space by explicitly constructing a semimetric $d$ on $\text{Vert}(X)$ in the following way. Let ${x, y \in \text{Vert}(X)}$ be any pair of vertices. If $x \neq y$, we set
\begin{equation} \label{semimetricequation}
d(x, y) := \begin{cases}
	\, \min  \{  i  \mid  \{x,y\} \in X_{i} \} \, \text{ if  } \,  \{x,y\} \in \text{Edge}(X),  \\
\infty \, \text{ otherwise},
	\end{cases}
\end{equation}
while if $x = y$, we set $d(x,y) := 0$. So the distance $d(x, y)$ between any two distinct vertices~$x, y$ is given by the smallest filtration step at which the edge $\{x, y\}$ enters into the filtration~$X_{\bullet}$, and it is assigned the value $\infty$ if $\{x, y\}$ is not an edge in $X$. Note that since $\mathbb{N}$ does not contain the number zero, we have $d(x,y)=0$ if and only if $x=y$. We remark that the distance function $d$ only defines a semimetric as it will in general not satisfy the triangle inequality.
With this understood, the \emph{Vietoris-Rips transformation} of $X$ is defined as the $\mathbb{N}$-filtered Vietoris-Rips complex
\begin{equation}
\widehat{\text{VR}}(X):=\text{VR}(\text{Vert}(X), d  )\notag
\end{equation}
with filtration $\widehat{\text{VR}}_{\bullet}(X)=(\text{VR}_{\, i} (\text{Vert}(X), d) )_{i \in \mathbb{N}}$. Our main result now states that the persistent homology in degree greater than zero of the $\mathbb{N}$-filtered flag complex $X$ can be computed in terms of the usual one-parameter persistent homology of its Vietoris-Rips transformation.
\begin{theorem}\label{mainthm} Let $X$ be a finite $\mathbb{N}$-filtered flag complex. Then in every positive degree $\ell > 0$, the Vietoris-Rips transformation induces an isomorphism
\begin{equation}
H_{\ell} (X ) \cong H_{\ell} \big (  \widehat{\text{\normalfont VR}}(X) \big) \label{VRIso}
\end{equation}
of persistence modules over $\mathbb{N}$.
\end{theorem}
Let us remark that the isomorphism \eqref{VRIso} will in general not hold in degree $\ell = 0$. This is because the filtration $\widehat{\text{VR}}_{\bullet}(X)$ is defined on the vertex set $\text{Vert}(X)$ of the whole complex, while the vertex set $\text{Vert}(X_i)$ can be a proper subset of $\text{Vert}(X)$.

\begin{proof}[Proof of the theorem]
By construction of the semimetric $d$ in \eqref{semimetricequation}, for every $i \in \mathbb{N}$ we have an identity
\[
X_{i}^{(1)} = \text{VR}_{\, i} (\text{Vert}(X), d)^{(1)} \setminus ( \text{Vert}(X) \setminus \text{Vert}(X_{i}) )
\]
of one-dimensional simplicial complexes, where on the right-hand side we need to remove all vertices in $X$ that are not contained in $X_{i}$.
Since both $X_{i}$ and $\text{VR}_{i} (\text{Vert}(X), d)$ are flag complexes, the above identity extends to an identity
\[
X_{i} = \text{VR}_{\, i}   (\text{Vert}(X), d  ) \setminus ( \text{Vert}(X) \setminus \text{Vert}(X_{i}) ).
\]
In particular, this identity means that the complexes $X_{i}$ and $\text{VR}_{\, i}  (\text{Vert}(X), d)$ consist of the same $k$-simplices for $k \geq 1$. Hence we immediately obtain the claimed isomorphism of persistence modules in all positive degrees $\ell > 0$.
\end{proof}

\subsection{Pathwise Persistence Barcodes}
\label{sec:PathwiseBarcodes}
A non-empty subset $I \subseteq \mathbb{N}$ is called an \emph{interval} if $r \leq s \leq t$ with $r,t \in I$ implies $s \in I$. We define a persistence module $\mathbb{F}_I$ over $\mathbb{N}$ as follows: \begin{equation}
\mathbb{F}_{I} (t) = \begin{cases}  \, \mathbb{F} \,\text{ if } t \in I,  \\
\, 0\, \text{ otherwise,}
\end{cases} \notag
\end{equation}
and $\mathbb{F}_I(s \leq t)$ is the identity map for all $s, t \in I$ with $s \leq t$ and the zero map otherwise. $\mathbb{F}_I$ is also called an \emph{interval module}.

Let $M$ be a persistence module over $\mathbb{N}$ of \textit{finite type}, i.e. $M$ is pointwise finite dimensional and there exists some $N \in \mathbb{N}$ such that $M(N \leq m)$ is an isomorphism of $\mathbb{F}$-vector spaces for all $m \in \mathbb{N}$ with $N \leq m$. The structure theorem for one-parameter persistence modules then states that $M$ admits a decomposition
\begin{equation}
M \cong \bigoplus_{j=1}^m\mathbb{F}_{I_j} \notag
\end{equation}
for a finite family $\mathcal{B}(M)=(I_1, \dots, I_m)$ of intervals which is uniquely determined up to the ordering of the intervals. This family $\mathcal{B}(M)$ is called the \emph{(persistence) barcode} of $M$. The intervals in the persistence barcode are also called \textit{bars}.

If $M=H_{\ell}(X)$ is the persistent homology in fixed degree $\ell \ge 0$ of some finite $\mathbb{N}$-filtered simplicial complex $X$, then $H_{\ell}(X)$ is of finite type and we denote its persistence barcode by $\mathcal{B}_{ \ell}(X)$. The persistence barcode $\mathcal{B}_{ \ell}(X)$ encodes the persistent homology of $X$ in degree $\ell$. The starting point of each bar in the persistence barcode corresponds to the birth of a homology feature, while its endpoint, if it exists, marks the death of the feature. The material about persistence barcodes summarized here is standard and can for example be found in \cite{Zomorodian2005} or \cite[\S 5.2]{Carlsson2005}.
\begin{figure}[t]
	\centering
	\includegraphics[width=1\textwidth]{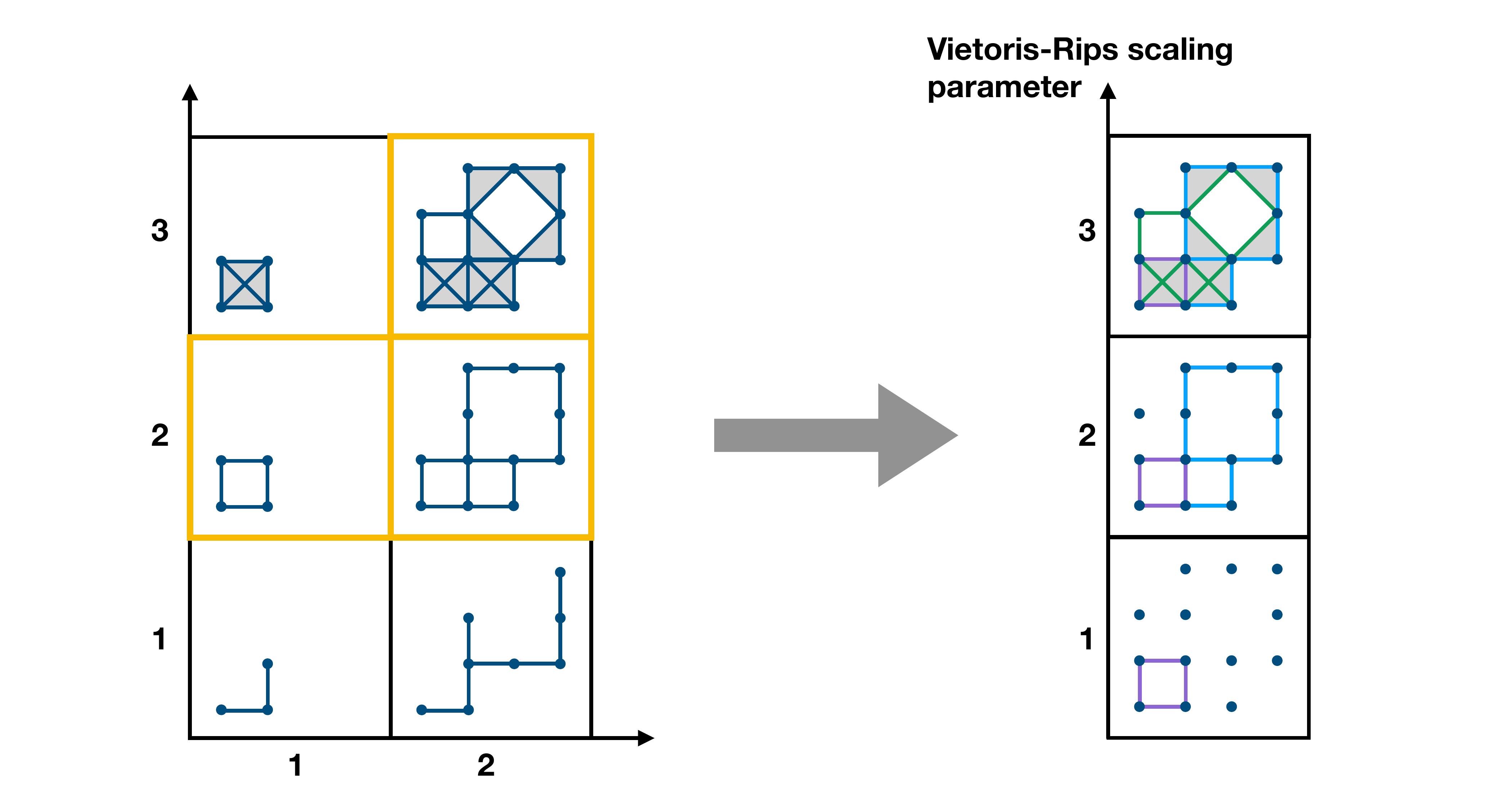}
\captionsetup{labelfont=bf}
\caption{\textbf{Example of a Vietoris-Rips transformation.}
On the left, we see a bi-filtered flag complex~$X$. The yellow squares mark the one-filtered subcomplex $X_{(1,2)} \subseteq X_{(2,2)} \subseteq X_{(2,3)}= X_{\nu}$ defined by the path $\nu=((1,2) \leq (2,2) \leq (2,3))$. The Vietoris-Rips transformation $\widehat{\text{VR}}(X_{\nu})$ of this subcomplex is the Vietoris-Rips complex shown on the right. The coloring of the edges indicates at which scale an edge enters into the Vietoris-Rips filtration.}
\label{fig:2}
\end{figure}

Let $P=(P,\leq)$ be a poset. A monotone sequence $\nu=(\nu_1 \leq \nu_2 \leq \nu_3 \leq \dots  )$ of elements in~$P$ is also called a \textit{(discrete) path} in $P$. The sequence $\nu$ \textit{stabilizes} if there exists some $m\in \mathbb{N}$ such that $\nu_i=\nu_m$ for all $i \geq m$. In this case, we use the notation $\nu=(\nu_1 \leq \dots \leq  \nu_m)$. Let us now consider a finite $P$-filtered flag complex~$X$.
Any path $\nu=(\nu_i)_{i \in \mathbb{N}}$ in $P$ gives rise to an $\mathbb{N}$-filtered subcomplex $X_{\nu}= \bigcup_{i \in \mathbb{N}} X_{\nu_i}$ of $X$ with filtration $(X_{\nu_i})_{i \in \mathbb{N}}$.
As an immediate consequence of \cref{mainthm}, in every positive degree~$\ell > 0$, the persistence barcode of $X_{\nu}$ may be computed in terms of a Vietoris-Rips persistence barcode of the Vietoris-Rips transformation as
\begin{equation}
\mathcal{B}_{  \ell}  (X_{\nu}) = \mathcal{B}_{ \ell}  \big (  \widehat{\text{\normalfont VR}} (X_{\nu}   ) \big ) . \notag
\end{equation}
We refer to $\mathcal{B}_{ \ell}  (X_{\nu})$ as the \emph{persistence barcode of $X$ along the path $\nu$}.
An instructive example of a Vietoris-Rips transformation is shown in Figure~\ref{fig:2}.

Our result demonstrates the usefulness of Vietoris-Rips transformations in practice. In fact, it reduces the computation of pathwise persistence barcodes in a finite multi-filtered flag complex to the computation of the persistence barcode of certain Vietoris-Rips filtrations.
In \cref{sec:MuRiT}, we will present the \texttt{MuRiT} algorithm, a software implementation of the Vietoris-Rips transformation for the efficient computation of pathwise persistence barcodes of multi-filtered flag complexes.

\subsection{Pathwise Persistence Barcodes and the Rank Invariant}
\label{subsec:RankInvariant}
Let $X$ be a finite $P$-filtered flag complex. Pathwise persistence barcodes are closely related to the rank invariant of $H_{\ell}(X)$ introduced by Carlsson and Zomorodian \cite{Carlsson2009multidimensionalpersistence}.
Let $P^2_{\Delta}=\{(v,w) \in P^2 \mid v \leq w\}  $.
Now the rank invariant of~$H_{\ell}(X)$ is given as the assignment
\begin{equation}
P^2_{\Delta} \to \mathbb{N}_{0}, \quad (v,w) \mapsto \text{rank}(H_{\ell}(X_v) \to  H_{\ell}(X_w)) \notag.
\end{equation}
It captures important persistence features of the multiparameter persistence module $H_{\ell}(X)$ and, as Carlsson and Zomorodian observerd, it is equivalent to the persistence barcode in the case of one-parameter persistence.
We can recover the rank invariant of $H_\ell(X)$ by computing the persistence barcode along the path $(\nu_1 \leq \nu_2 )$ for every pair $(\nu_1,\nu_2) \in P^2_{\Delta}$.

\subsection{The \texttt{MuRiT} Algorithm for Multi-Filtered Flag Complexes}\label{sec:MuRiT}
Based on our theoretical considerations in the previous subsections, we now introduce the \texttt{MuRiT} algorithm, displayed in \cref{alg:murit}. It is designed for the computation of pathwise persistence barcodes in positive homology degree of finite multi-filtered flag complexes via Vietoris-Rips transformations. The \texttt{MuRiT} algorithm firstly efficiently computes the Vietoris-Rips transformation, and secondly uses \texttt{Ripser} to compute the persistence barcodes of this Vietoris-Rips transformation.

Our setup for \texttt{MuRiT} will be a finite $P$-filtered flag complex $X$ with filtration $X_{\bullet}=(X_{p})_{p \in P}$ for some finite $n$-dimensional subposet $P \subseteq \mathbb{R}^n$.
This ensures that \texttt{MuRiT} will be applicable to a large class of real-world data, like for example time series data.
Instead of encoding the full complex~$X$, it will be sufficient to specify a finite \emph{edge entry annotation list} \begin{equation}
L: \text{Edge}(X) \to \mathcal{P}(P) \notag
\end{equation} which takes values in the power set $\mathcal{P}(P)$ of $P$ and records, for every edge $\{x,y\} \in \text{Edge}(X)$, the minimal filtration steps~$p\in P$ at which this edge enters into the filtration:
\begin{align}
	L( \{x,y\} ) := \min  \{ p \in P \mid \{x,y\} \in X_p  \} \subseteq P \notag
\end{align}
Recall at this point that minima of subsets of posets are in general not unique. Lastly, we need to specify a path $\nu=(\nu_1 \leq \dots \leq \nu_m)$ in $P$, which defines the one-filtered subcomplex $X_\nu \subseteq X$ we would like to analyze.

From this input data, \texttt{MuRiT} first computes the lower triangular distance matrix $D$ of the following semimetric $d$ on the vertex set $\text{Vert}(X)$: the restriction of $d$ to the vertex set $\text{Vert}(X_{\nu})$ coincides with the semimetric \eqref{semimetricequation} associated with the Vietoris-Rips transformation $\widehat{\text{\normalfont VR}}(X_{\nu})$, while for any pair $x,y \in \text{Vert}(X) \setminus \text{Vert}(X_{\nu})$, we set $d(x,y):=\infty$ if $x\neq y$ and $d(x,y)=0$ if $x=y$.
This makes the algorithm more user friendly---the user only has to encode the complex~$X$ once by specifying the annotation list $L$.
After that, they only have to define the paths in $P$ along which they want to compute persistence barcodes.

To fix the notation, we denote the vertices in $X$ by $\text{Vert}(X)=\{x_1, \dots, x_N\}$. Then $D$ is given by $D_{ij}=d(x_i,x_j)$ with $i>j$.
In order to determine the distance~$D_{ij}$, \texttt{MuRiT} calculates the unique minimal intersection of the upper set of $L(\{x_{i},x_{j}\})$ with the given path $\nu$ in $P$.
Recall that for a subset $Q \subseteq P$, the upper set of $Q$ in $P$ is defined as the set of all $p\in P$ such that $q \leq p$ for some~$q \in Q$.
In a second step, \texttt{MuRiT} passes the distance matrix $D$ to \texttt{Ripser} for the actual computation of the persistence barcodes $\mathcal{B}_{\bullet}(X_{\nu}) := (\mathcal{B}_{\ell }(X_{\nu}))_{\ell \geq 1}$ in positive homology degree. Note at this point that the distance matrix $D$ will in general not satisfy the triangle inequality. 
But this is not a problem as \texttt{Ripser} can handle distance matrices of semimetric spaces and in particular does not require the matrix $D$ to satisfy the triangle inequality \cite{notriangle}.

\begin{algorithm}[h]
\caption{\texttt{MuRiT} algorithm}\label{alg:murit}

\begin{algorithmic}[1]
\Require
\Statex
Vertex Set $\text{Vert}(X)=\{x_1, \dots, x_N\}$
\Statex Edge Entry Annotation List $L: \text{Edge}(X) \to \mathcal{P}(P)$
\Statex Path $\nu=(\nu_1 \leq \dots \leq \nu_m)$ in $P$
\vspace{0.5em}
\Ensure
\Statex Persistence barcodes $\mathcal{B}_{\bullet}(X_{\nu})$
\vspace{1em}
\ParallelFor{every $x_i,x_j$ in $\text{Vert}(X)$ with $i > j$} \hspace{1em} \tikzmark{right}
\State $D(x_i,x_j) := \infty$
	\For{every step $\nu_k$ in the path $\nu$}
	\For{every filtration step $p$ in the edge entry annotation $L(\{x_i,x_j\})$}
			\If{$p \leq \nu_k$} \Comment{check if the edge $\{x_i,x_j\}$ is contained in $X_{\nu_{k}}$}
				\State $D(x_i,x_j) = k$
				\State \Goto \texttt{~bottom}
\EndIf
		\EndFor
	\EndFor

	\State \texttt{bottom}
\EndParallelFor

\vspace{0.2em}

\State $\mathcal{B}_{\bullet}(X_{\nu}) := $ \Call{Ripser}{$D$}
\State \Return{$\mathcal{B}_{\bullet}(X_{\nu})$}
\end{algorithmic}
\end{algorithm}

A common way in which multi-filtered flag complexes naturally arise in applications is in the Vietoris-Rips persistence analysis of multi-filtered metric datasets that come with extra structure given by additional parameters. A typical example of this is time series data in the evolution of the coronavirus (see \cref{sec:application}). To formalize this, let $(S,h)$ be a finite semimetric space equipped with a filtration $S_{\bullet}=(S_{{t}})_{{t} \in T}$ for some finite $n$-dimensional subposet $T \subseteq \mathbb{R}^n$. For example, in the case of time series data we could have a totally ordered subset $T = \{ t_{1} \leq \ldots \leq t_{m} \} \subseteq \mathbb{R}$ that specifies the time steps.
We may then consider the Vietoris-Rips complex $\text{VR}_{\,r}(S_{{t}},h)$ for each filtration step~${t} \in T$ and $r \in h(S)$, where $h(S)$ denotes the set of pairwise distances $h(x,y)$ of elements $x,y \in S$.
This gives rise to a multi-filtered flag complex in the following way.
Consider the poset
\[
P := T \times h(S) \subseteq \mathbb{R}^{n+1}.
\]
Then $\text{VR}(S,h)$ naturally becomes a $P$-filtered flag complex with filtration
\begin{equation}
\text{VR}_{\bullet}(S_{\bullet},h)=   ( \text{VR}_{\, r}(S_{{t}},h)  )_{({t},r) \in P} \notag
\end{equation}
As a result, we may use \texttt{MuRiT} to efficiently investigate the multi-parameter persistent homology of the multi-filtered flag complex $\text{VR}(S,h)$ by computing pathwise persistence barcodes via Vietoris-Rips transformations.

In order to be able to run the \texttt{MuRiT} \cref{alg:murit}, we first need to prepare the input data, which will be done with \cref{alg:metric+filtration-to-edge-annotation}.
We denote the points in the dataset by $S=\{x_1, \dots, x_N\}$. Moreover, we denote by $H$ the lower triangular distance matrix of the semimetric space $(S,h)$ given by $H_{ij}=h(x_i,x_j)$ with ${i>j}$.
\cref{alg:metric+filtration-to-edge-annotation} takes as input this lower triangular distance matrix $H$, a path $\nu$ in $P$ of the form ${\nu=(({t}_1, r_1) \leq \dots \leq ({t}_m, r_m))}$, and a \emph{point entry annotation list} $K: S \to \mathcal{P}(T)$ where for every point $x \in S$
\begin{equation}
	K(x) := \min  \{ {t} \in  T  \mid x \in S_{{t}}  \} \subseteq T \notag.
\end{equation}
The output of \cref{alg:metric+filtration-to-edge-annotation} is the edge entry annotation list $L_{\nu}: \text{Edge}(\text{VR}(S,h)_{\nu}) \to \mathcal{P}(P)$ of the one-filtered subcomplex $\text{VR}(S,h)_{\nu}\subseteq \text{VR}(S,h)$ given by
\[
L_{\nu}( \{x,y\} ) = \min  \{ ({t}_i,r_i) \in \nu \mid \{x,y\} \in \text{VR}_{r_i}(S_{{t}_i},h)  \} \subseteq P.
\]
We may then pass as input for \cref{alg:murit} the one-filtered subcomplex $\text{VR}(S,h)_{\nu}\subseteq \text{VR}(S,h)$ with vertex set $S_{{t}_m} \subseteq S$, together with the annotation list $L_{\nu}$. 

\begin{algorithm}[h]
\caption{Preparation of Multi-Filtered Data for \texttt{MuRiT}}\label{alg:metric+filtration-to-edge-annotation}
\begin{algorithmic}[1]
\Require
\Statex Lower Triangular Distance Matrix $H$
\Statex Point Entry Annotation List $K: S \to \mathcal{P}(T)$
\Statex Path $\nu=(({t}_1, r_1) \leq \dots \leq ({t}_m, r_m))$ in $P$
\vspace{0.5em}
\Ensure
\Statex Edge Entry Annotation List $L_{\nu}: \text{Edge}(\text{VR}(S,h)_{\nu}) \to \mathcal{P}(P)$
\vspace{1em}

\Function{GetPointOfEntry}{List $K(x)$, path $\nu$}
\For{every step $({t}_i, r_i)$ in the path $\nu$}
	\For{every filtration step $u$ in the point entry annotation $K(x)$}
		\If{${u} \leq {t}_i$} \Comment{check if the point $x$ is contained in $S_{t_{i}}$}
			\State \Return $({t}_i, r_i)$
		\EndIf
	\EndFor
\EndFor
\State \Return $\infty$
\EndFunction

\Statex

\Procedure{GetEdgeEntryAnnotation}{}
\State Set $L_\nu$ to empty list
\ParallelFor{every pair of points $x_{i},x_{j}$ in $S$ with $i > j$}
	\State $a$ := \Call{GetPointOfEntry}{$K(x_i)$, $\nu$}
	\State $b$ := \Call{GetPointOfEntry}{$K(x_j)$, $\nu$}
	\State Append $(\max(a, b), H_{ij})$ to $L_{\nu}$
\EndParallelFor
\EndProcedure
\end{algorithmic}
\end{algorithm}
A parallelized implementation of the \texttt{MuRiT} algorithm, is available at \url{https://github.com/tdalife/murit}.
This implementation combines Algorithms \ref{alg:murit} and \ref{alg:metric+filtration-to-edge-annotation}, and is optimized for the Vietoris-Rips persistence analysis of multi-filtered datasets.

\section{An Application to Viral Evolution}
\label{sec:application}

\subsection{Topology of Evolutionary Processes}
A particularly useful application of the MuRiT algorithm is in the surveillance of pathogen evolution proposed in \cite{topologyidentifies2022}.
While the evolution of an organism is usually modeled according to the paradigm of phylogenetic trees, there are various biological phenomena that are incompatible with this approach.
In the example of viral evolution, an exchange of genomic information between distinct lineages can happen through recombination or reassortment, which is known to be a key factor in rapid host adaptation for many viruses \cite{Chan2013}.
These instances of reticulate evolution can be viewed as a deviation from a trivial tree topology and obstruct the existence of a single phylogeny, as different parts of the genome might admit conflicting evolutionary histories.

In this context, we consider a finite set $S$ of genome sequences, which we will think of as finite words {$x=(x_{1}, x_{2}, \ldots, x_{N})$} of uniform length $N$ in the alphabet $\{\text{A},\text{C},\text{T},\text{G}\}$ corresponding to the four nucleotides, and endow it with a semimetric measuring the genetic distance between any pair of sequences.
A standard choice of semimetric is given by the Hamming distance
\[
h: S \times S \to  \mathbb{N}\qquad (x,y)\mapsto h(x,y):= \#\{i\,|\,x_i\neq y_i\},
\]
which counts the number of genomic positions at which two sequences differ.
As put forward in \cite{Chan2013}, persistent homology provides a fast and effective method to extract patterns of non-trivial topology from the genomic data set $S$.
In this approach, evolutionary relationships at all scales are comprehensively encoded in the Vietoris-Rips filtration
\begin{equation}\text{VR}_{\, 0} (S,h) \subseteq \text{VR}_{\, 1} (S,h) \subseteq \text{VR}_{\,  2} (S,h) \subseteq \dots, \notag
\end{equation}
while information on both the tree-structure and reticulate events is captured in the persistence barcode.

\subsection{Topological Recurrence Analysis}\label{sec:Topological recurrence analysis}

In addition to the exchange of genetic material, a further source of non-trivial topology in the phylogeny is convergent evolution.
The corresponding reticulate events are called \textit{homoplasies} and mean the independent acquisition of a specific mutation in different lineages.
If the sampling of data is sufficiently dense, this typically gives rise to a persistent homology class at the smallest scale, admitting a cycle representative for which all edges correspond to \textit{single nucleotide variations} (SNV) with edge length equal to~$1$.
Recall at this point that a single nucleotide variation means a mutation of a single nucleotide at precisely one position in the genome.
Such cycles admit a convenient description in terms of their location in the barcode.
The representatives of persistent homology classes corresponding to bars in the barcode $\mathcal{B}_{ 1} (\text{VR}(S,h))$ that are born in the first filtration step $ \text{VR}_{\, 1}(S,h)$ will be called \textit{SNV cycles}.

In \cite{topologyidentifies2022}, SNV cycles are used to define a measure of convergent evolution based on one-parameter persistent homology.
Let us briefly outline the main ideas of this approach.
The barcode associated with the one-dimensional persistent homology $H_1(\text{VR}(S,h))$ is computed with a custom version of \texttt{Ripser} \cite{bauer2021ripser} that also has the capability to find cycle representatives.
This method of \textit{exhaustive reduction} (see \cite{Zomorodian2005, Edelsbrunner2019}) aims to produce cycles that tightly fit the data, by systematically replacing the longest edge of a given cycle with shorter edges of suitable nearby cycles.
By definition of SNV cycles, for their extraction it suffices to consider persistent homology only on the smallest scale, which amounts to running \texttt{Ripser} with scale parameter threshold set to~$2$, leading to a substantial speedup of the computation.
Thanks to the specific properties of the Hamming geometry of the semimetric space $(S,h)$ in combination with the tree-like structure of the phylogeny, \texttt{Ripser} is able to process hundreds of thousands of data points in the particular case of SARS-CoV-2 evolution \cite{topologyidentifies2022, gromovhyper}.

For each mutation in the viral genome, its \textit{topological recurrence index} (tRI) is defined to be the total number of SNV cycles containing an edge that gives rise to the given mutation.
Here we use the fact that each edge in an SNV cycle has length $1$ and hence corresponds to a uniquely determined mutation in the genome.
The topological recurrence index provides a lower bound for the number of independent occurrences of a mutation in the evolution of an organism and is therefore a measure for convergent evolution.
As was demonstrated in \cite{topologyidentifies2022}, there is an abundance of SNV cycles in the SARS-CoV-2 genomic dataset and hence the above definition of the topological recurrence index, which only relies on SNV cycles, will already lead to statistical significant signals.
It was moreover shown that the topological recurrence index can serve as an early indicator for emerging adaptive mutations in the evolution of the coronavirus.

In the case of ongoing genomic surveillance, the dataset admits a natural filtration by sampling time $S_1\subseteq \cdots S_{m-1}\subseteq S_m$, where $S_t$ denotes the set of all viral genomes sampled up until time step $t$.
For every $t \in \{1, \dots, m\}$, we denote by $\textbf{SNV}_t$ the full set of SNV cycle representatives in time step $t$ we get from the persistence analysis of $(S_{t},h)$ with \texttt{Ripser}.
Note that, as a feature of Hamming geometry, the homology classes associated to SNV cycles are in a certain sense stable with respect to adding points in a time series:
Inclusion of new data points might lead to a splitting of the cycle, but the resulting homology class is typically non-zero. It can happen that the homology class is destroyed by adding data points only when gaps or insertions in the genome alignment are involved in the SNV cycle.

In the study of SARS-CoV-2 evolution in \cite{topologyidentifies2022}, the barcode $\mathcal{B}_{ 1}( \text{VR}(S_t,h))$ for each time step~$t$ is computed separately.
Tracking SNV cycles over time yields information about the adaptation process of the pathogen. However, this can be troubled by the following two issues.

\begin{enumerate}
\item[(1)] Computing persistence homology at each time step separately can be time consuming and computationally expensive if the filtration by time of the dataset is large (like for example in a time series analysis over one year on a daily basis).
\item[(2)] SNV cycle representatives of the same homology class at different time steps will in general not be compatible with each other, which can lead to noise in the topological recurrence analysis. More precisely, if $\omega  \in \textbf{SNV}_t$ and if its image under the canonical morphism
\begin{equation}
H_1( \text{VR}_{\, 1}(S_t,h)) \to H_1( \text{VR}_{\, 1}(S_{t+1},h)) \notag
\end{equation}
induced by the inclusion $S_t \subseteq S_{t+1}$ is non-zero, then it may still happen that $\omega \not \in \textbf{SNV}_{t+1}$.
\end{enumerate}

\begin{figure}[H]
\centering
	\includegraphics[width=1\textwidth]{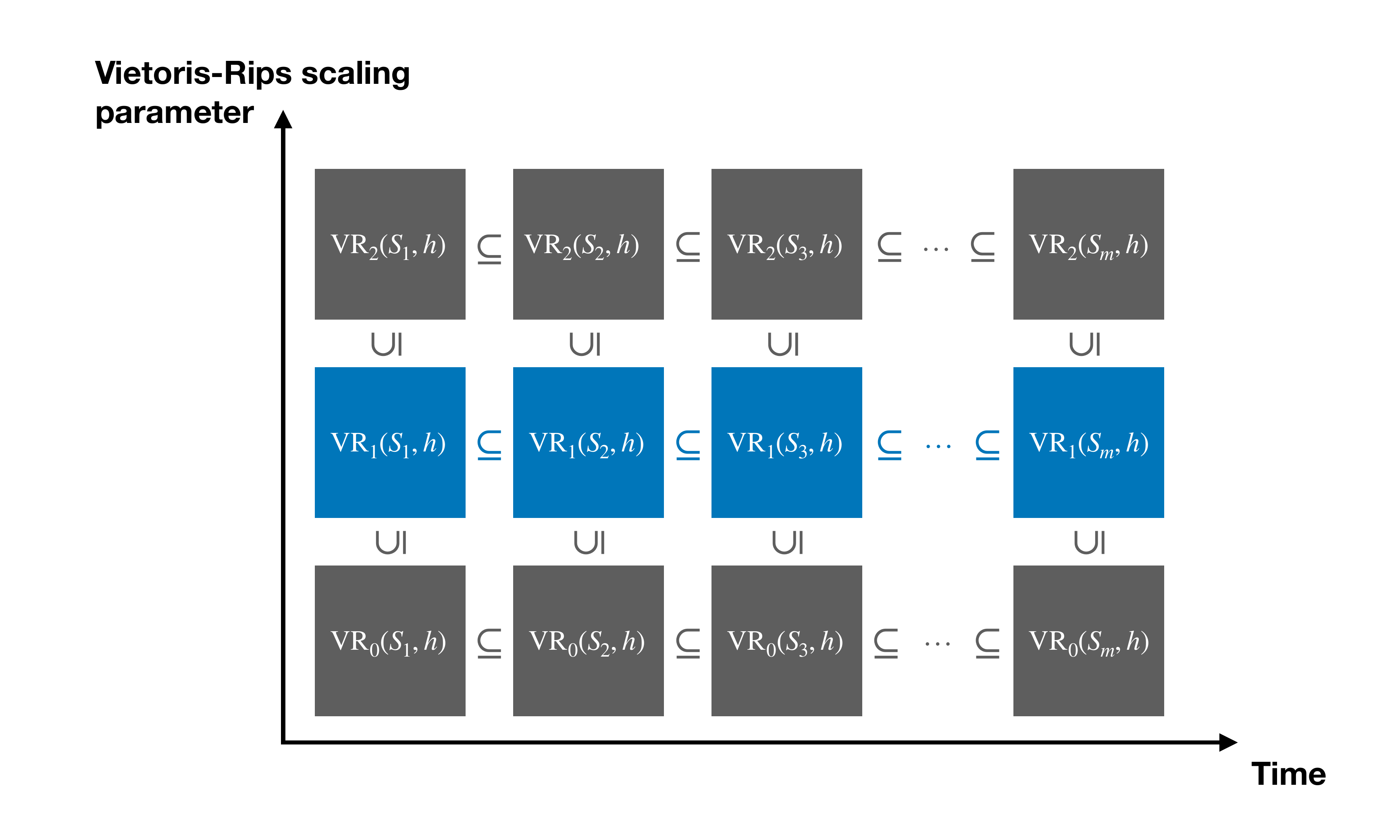}
\captionsetup{labelfont=bf}
\caption{\textbf{Pathwise subcomplexes and SNV cycles in viral evolution.}
The blue tiles mark the horizontal subcomplex $\text{VR}(S,h)_{\nu}$ of the $P$-filtered Vietoris-Rips complex $\text{VR}_{\bullet}(S_{\bullet},h)$ corresponding to the path $\nu=((1,1) \leq \dots \leq (m,1))$ in $P = \{1, \dots, m \} \times h(S)$. This subcomplex keeps track of the formation of SNV cycles in the time-filtered dataset $S_1\subseteq S_{2} \subseteq \cdots \subseteq S_m = S$ of viral gene sequences equipped with the Hamming distance $h$.}
\label{fig:3}
\end{figure}

Both of these problems are resolved by the MuRiT algorithm in the following way. As explained in \cref{sec:MuRiT}, we define the poset
\begin{equation}
P :=\underbrace{ \{1, \dots, m \}}_{\text{time steps}} \times \underbrace{h(S)}_{\text{distances}} \notag.
\end{equation}
Then the Vietoris-Rips complex $X=\text{VR}(S,h)$ naturally becomes a $P$-filtered simplicial complex. As for the definition of the topological recurrence index it is only relevant to determine the time of birth of each SNV cycle, it suffices to compute the persistence barcode in homology degree one of the subcomplex $\text{VR}(S,h)_{\nu} \subseteq \text{VR}(S,h)$ determined by the path $\nu=((1,1) \leq \dots \leq (m,1))$ in $P$ (see Figure~\ref{fig:3}).
For this, we use \texttt{MuRiT} to compute the persistence barcode in homology degree one of the Vietoris-Rips transformation of $\text{VR}(S,h)_{\nu}$.
At each time step $t$, \texttt{Ripser} extracts a full set of SNV cycles such that the corresponding homology classes correspond to bars $I  \in \mathcal{B}_1 \big( \widehat{\text{VR}}(\text{VR}(S,h)_{\nu})  \big )$ with $t \in I$.
In this way, we achieve compatibility of SNV cycles across all filtration steps.
From a computational perspective, this means a great gain of efficiency as the persistence barcode of $H_{1}(\text{VR}(S,h)_{\nu})$, which resolves all time steps, is computed with a single run of \texttt{Ripser}.

\subsection{The CoVtRec Pipeline}
\label{sec:covtrec}
The \texttt{MuRiT} algorithm is part of our \texttt{CoVtRec} pipeline for the early warning and surveillance of recurrent mutations in the evolution of the coronavirus {SARS-CoV-2} in the current COVID-19 pandemic \cite{covtrec}. Regular reports containing analyses on the basis of SARS-CoV-2 genomic data shared via GISAID, the global data science initiative \cite{Shu2017, Khare2021}, are available at \url{https://tdalife.github.io/covtrec}.
Recall that recurrent mutations are potentially adaptive in the sense that they could confer some fitness advantage to the virus, such as immune evasion or higher infectivity.
In the current phase of the pandemic, the early identification of potentially adaptive mutations is of paramount importance as the virus is constantly developing new variants by mutating its genome. For more details on the biological aspects of the topological recurrence analysis of SARS-CoV-2 genomic data see \cite{topologyidentifies2022}.

\begin{figure}[h!]
\captionsetup{format=plain, labelfont={bf},labelformat={default},labelsep=period,name={Figure}}
\centering
\includegraphics[scale=0.32]{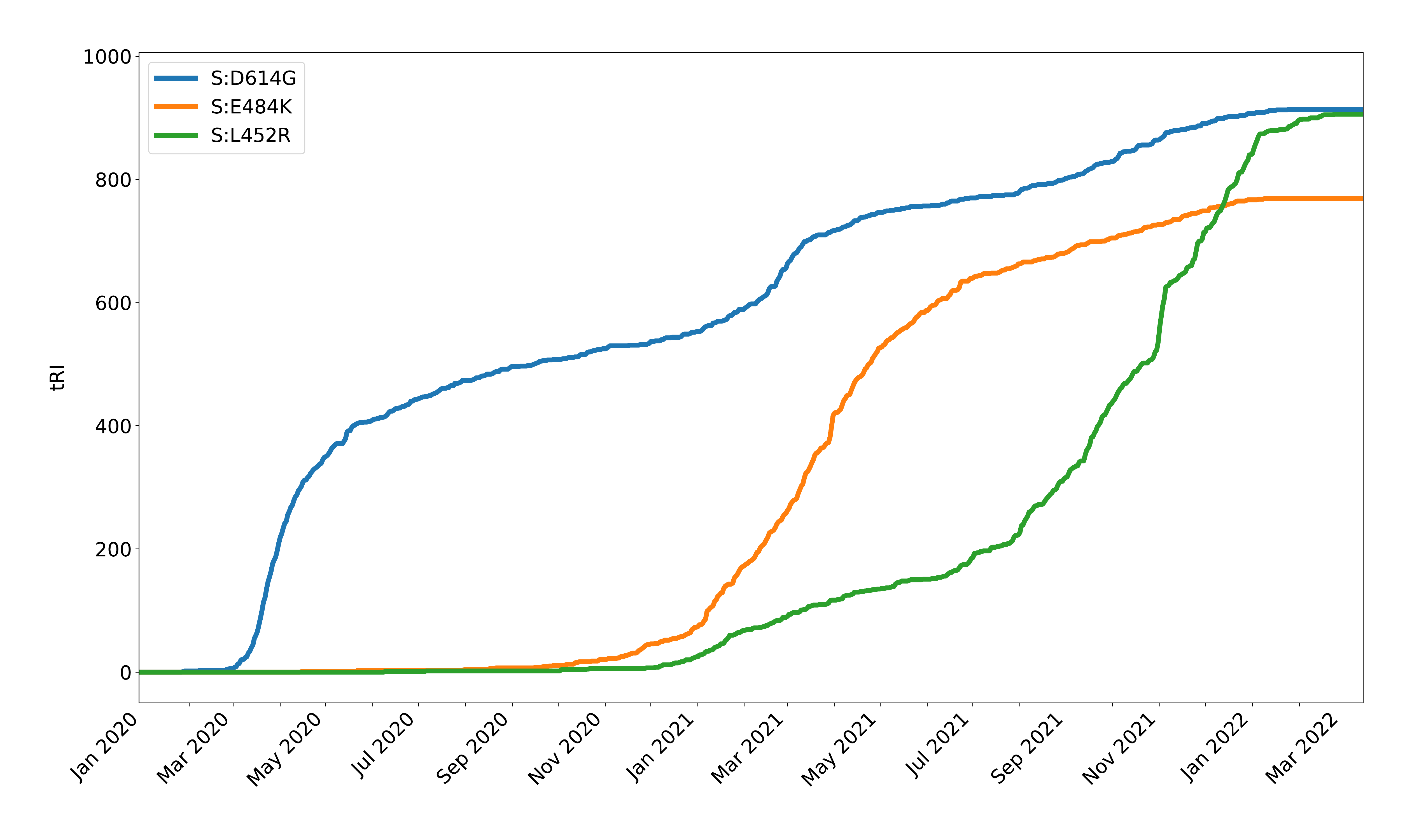}
\caption{\textbf{Surveillance of the convergent evolution of the coronavirus SARS-CoV-2.}
The diagram shows the time plots at daily resolution for the topological recurrence index (tRI) of the adaptive SARS-CoV-2 Spike gene mutations \texttt{D614G}, \texttt{E484K} and \texttt{L452R} over a period of~27 months, from the beginning of the pandemic in late December 2019 until 15 March 2022.
The \texttt{MuRiT} algorithm enables the efficient topological analysis of hundreds of thousands of data points over time by leveraging the natural stratification by time of the coronavirus gene sequences dataset.
While the mutation \texttt{D614G} is currently observed in essentially all virus samples, \texttt{E484K} occurred in the Alpha and Beta variants, and \texttt{L452R} has more recently been observed in the Delta and Omicron BA.5 variants.}
\label{fig:plot}
\end{figure}

The \texttt{CoVtRec} pipeline generates time series analysis charts for the topological recurrence index (tRI) at daily resolution by leveraging the natural stratification by time of SARS-CoV-2 genomic data.
Thanks to highly optimized algorithms that take advantage of the tree-like structure of the data \cite{gromovhyper}, \texttt{CoVtRec} can process very large SARS-CoV-2 genomic datasets and easily scales to hundreds of thousands of distinct genomes.
This demonstrates the efficiency and usefulness of the \texttt{MuRiT} algorithm in practice.
To give a concrete example, we analyzed topological signals for the ongoing convergent evolution for three prominent mutations of the SARS-CoV-2 genome from the beginning of the pandemic in December~2019 until \Date.
To that end, we performed a topological recurrence analysis for a curated alignment of \NumberMSACleanSequences high-quality SARS-CoV-2 Spike gene sequences shared via GISAID.
The analysis was restricted to the Spike gene, a part of the genome that determines the structure of the Spike protein on the surface of the virus and therefore plays an essential role in immune evasion and infectivity.
Our algorithm performed the topological analysis of \NumberMSADistinctSequences distinct Spike gene sequences in less than a day on a machine with Intel Xeon Gold 6230R processors and 52 kernels.

We analyzed topological signals of convergence for the Spike gene mutations \texttt{D614G}, \texttt{E484K} and \texttt{L452R} (see Figure \ref{fig:plot}). All of these mutations exhibit a topological signal of convergence, with the topological recurrence index (tRI) rising over the course of the pandemic. We conclude that they are potentially adaptive. In fact, there is by now experimental evidence that the mutation \texttt{D614G} increases transmissibility \cite{Korber2020, Li2020} and in vitro infectiousness \cite{Plante2021, Hou2020, Yurkovetskiy2020}, and the mutations \texttt{E484K} and \texttt{L452R} enable the virus to evade immune protection \cite{Greaney2021, Liu2021a}. While the mutation \texttt{D614G} superseded the wild type already at the beginning of the pandemic and is currently observed in essentially all virus samples, \texttt{E484K} occurred in the Alpha and Beta variants, and \texttt{L452R} has more recently been observed in the Delta and Omicron BA.5 variants \cite{WHO2021}. A more detailed discussion of biological implications of the presence of topological signals for mutations in the evolution of the coronavirus is available in \cite{topologyidentifies2022}.

\subsection*{Data Availability and Data Preparation} All SARS-CoV-2 genome data used in this work are available from the GISAID EpiCov Database \cite{Shu2017, Khare2021} and are accessible online at \episetDOI.
Our analysis of SARS-CoV-2 genome data is based on the alignment \Alignment {} downloaded from the GISAID EpiCoV Database \cite{Shu2017, Khare2021} on \DownloadDate. This alignment comprises \NumberMSASequences SARS-CoV-2 whole genome sequences that have been aligned to the reference sequence Wuhan/WIV04 with GISAID accession number \texttt{EPI\_ISL\_402124} using \texttt{MAFFT (Version 7)} \cite{Katoh2002}.
Sequences in this alignment were truncated to the Spike gene (reference site positions 21,563 to 25,384), and subsequently sequences containing any characters other than A, C, T, G or gaps or insertions represented by - were removed. This resulted in an alignment comprising \NumberMSACleanSequences complete SARS-CoV-2 Spike genes of length \LengthMSACleanSequences.

\bibliographystyle{aomalpha}
\bibliography{article}

\end{document}